\documentclass[a4paper,11pt]{article}
\usepackage{amsmath}
\usepackage{amssymb}
\usepackage{amsthm}
\usepackage{mathrsfs}
\usepackage[russian]{babel}
\usepackage{graphics}

\newcommand{\al}{\alpha}
\newcommand{\be}{\beta}
\newcommand{\RR}{\mathbb R}

\newcommand{\ZZ}{\mathbb Z}
\newcommand{\EE}{\mathbb E}
\renewcommand{\P}{{\cal P}}
\newcommand{\K}{\tilde K}
\newcommand{\G}{\mathscr D}
\newcommand{\D}{\mathscr H}
\renewcommand{\H}{\mathscr K}

\newcommand{\ind}{{\rm ind\,}}
\newcommand{\eps}{\varepsilon}
\newcommand{\ups}{\upsilon}
\newcommand{\grad}{{\rm grad\,}}
\newcommand{\rank}{{\rm rank\,}}
\newcommand{\id}{{\rm id}}
\newcommand{\isot}{{\rm isot}}
\newcommand{\abs}{{\rm abs}}

\newcommand{\glue}{{\rm glue}}
\newcommand{\col}{\colon\,}
\renewcommand{\d}{\partial}
\newcommand{\Proof}{\noindent{\bf Доказательство. }}
\newcommand{\PProof}[1]{\noindent{\bf Доказательство {#1}. }}
\newcommand{\la}{\langle}
\newcommand{\ra}{\rangle}
\newcommand{\<}{\langle\!\langle}
\renewcommand{\>}{\rangle\!\rangle}
\renewcommand{\tilde}{\widetilde}
\renewcommand{\hat}{\widehat}

\theoremstyle{plain}
\newtheorem{Thm}{Теорема}
\newtheorem*{Sta}{Утверждение}

\newtheorem*{Cor}{Следствие}
{\theoremstyle{definition}
                             
                             \newtheorem{Def}{Определение}
                             \newtheorem{Not}{Обозначение}

                             \newtheorem*{Exs}{Примеры}
}

\catcode`\@=11
\newcommand{\smalllineskip}{\baselineskip=10pt}
\newcommand{\fcaption}[1]{
        \refstepcounter{figure}
        \setbox\@tempboxa = \hbox{\footnotesize {\bf Рис.~\thefigure.} #1}
        \ifdim \wd\@tempboxa > 5in
           {\begin{center}
        \parbox{5in}{\footnotesize\smalllineskip {\bf Рис.~\thefigure.} #1}
            \end{center}}
        \else
             {\begin{center}
             {\footnotesize {\bf Рис.~\thefigure.} #1}
              \end{center}}
        \fi}

\title{Connected components of spaces of Morse functions \\
with fixed critical points}
\author{Elena A.\ Kudryavtseva}
\date{}

\begin{document}

\maketitle

{\small {\bf Abstract.} Let $M$ be a smooth closed orientable surface
and $F=F_{p,q,r}$ be the space of Morse functions on $M$ having
exactly $p$ critical points of local minima, $q\ge1$ saddle critical
points, and $r$ critical points of local maxima, moreover all the
points are fixed. Let $F_f$ be the connected component of a function
$f\in F$ in $F$. By means of the winding number introduced by
Reinhart (1960), a surjection $\pi_0(F)\to\ZZ^{p+r-1}$ is
constructed. In particular, $|\pi_0(F)|=\infty$, and the Dehn twist
about the boundary of any disk containing exactly two critical
points, exactly one of which is a saddle point, does not preserve
$F_f$. Let $\G$ be the group of orientation preserving
diffeomorphisms of $M$ leaving fixed the critical points, $\G^0$ be
the connected component of $\id_M$ in $\G$, and $\G_f\subset\G$ the
set of diffeomorphisms preserving $F_f$. Let $\D_f$ be the subgroup
of $\G_f$ generated by $\G^0$ and all diffeomorphisms $h\in\G$ which
preserve some functions $f_1\in F_f$, and let $\D_f^\abs$ be its
subgroup generated by $\G^0$ and the Dehn twists about the components of
level curves of functions $f_1\in F_f$. We prove that
$\D_f^\abs\subsetneq\G_f$ if $q\ge2$, and construct an epimorphism
$\G_f/\D_f^\abs\to\ZZ_2^{q-1}$, by means of the winding number. A
finite polyhedral complex $K=K_{p,q,r}$ associated to the space $F$
is defined. An epimorphism $\mu\col\pi_1(K)\to\G_f/\D_f$ and finite
generating sets for the groups $\G_f/\G^0$ and $\G_f/\D_f$ in terms
of the 2-skeleton of the complex $K$ are constructed.


{\it MSC-class: } 58E05, 57M50, 58K65, 46M18
}

 \begin{center} \bf {\LARGE
Связные компоненты пространств функций Морса \\
с фиксированными критическими точками} \medskip\\
Елена А.\ Кудрявцева
 \end{center}

\begin{abstract}
Пусть $M$ -- гладкая замкнутая ориентируемая поверхность и
$F=F_{p,q,r}$ -- пространство функций Морса на $M$, имеющих ровно $p$
критических точек локальных минимумов, $q\ge1$ седловых критических
точек и $r$ точек локальных максимумов, причем эти точки фиксированы.
Пусть $F_f$ -- компонента связности функции $f\in F$ в $F$. С помощью
числа вращения, введенного Рейнхартом (1960), построена сюръекция
$\pi_0(F)\to\ZZ^{p+r-1}$. В частности, $|\pi_0(F)|=\infty$, и
скру\-чи\-ва\-ние Дэна вокруг границы любого диска, содержащего ровно
две критические точки, из которых ровно одна седловая, не сохраняет
$F_f$. Пусть $\G$ -- группа сохраняющих ориентацию диффеоморфизмов
$M$, оставляющих неподвижными критические точки, $\G^0$ -- компонента
связности $\id_M$ в $\G$, $\G_f\subset\G$ -- множество
диффеоморфизмов, сохраняющих $F_f$. Пусть $\D_f$ -- подгруппа $\G_f$,
порожденная $\G^0$ и всеми диффеоморфизмами $h\in\G$, сохраняющими
какие-либо функции $f_1\in F_f$, и пусть $\D_f^\abs$ -- ее подгруппа,
порожденная $\G^0$ и скручиваниями Дэна вокруг компонент линий уровня
функций $f_1\in F_f$. С помощью числа вращения доказано, что
$\D_f^\abs\subsetneq\G_f$ при $q\ge2$, и построен эпиморфизм
$\G_f/\D_f^\abs\to\ZZ_2^{q-1}$. Определен конечный полиэдральный
комплекс $K=K_{p,q,r}$, ассоциированный с пространством $F$.
Построены эпиморфизм $\mu\col\pi_1(K)\to\G_f/\D_f$ и конечные
множества порождающих элементов групп $\G_f/\G^0$ и $\G_f/\D_f$ в
терминах 2-остова комплекса $K$.

\medskip
{\it Ключевые слова: } функции Морса на поверхности, эквивалентные и
изотопные функции, число вращения, скручивание Дэна, допустимый
диффеоморфизм, полиэдральный комплекс.

{\it УДК } 515.164.174, 515.122.55
\end{abstract}


{\bf 1. Введение.} Пусть $M=M^2_g$ -- гладкая замкнутая ориентируемая
поверхность и $F=F_{p,q,r}$ -- пространство функций Морса на $M$,
имеющих ровно $q\ge1$ седловых критических точек $x_1,\dots,x_q$, $p$
критических точек $x_{q+1},\dots,x_{p+q}$ локальных минимумов и $r$
точек $x_{p+q+1},\dots,x_{p+q+r}$ локальных максимумов, причем эти
точки фиксированы. Возникает задача: описать гомотопический тип
про\-стран\-ства $F$ (в $C^\infty$-топологии) и, в частности,
множество $\pi_0(F)$ его связных компонент. С помощью {\it числа
вращения}, введенного Б.~Рейнхартом~\cite {R}, мы строим сюръекцию
$\pi_0(F)\to\ZZ^{p+r-1}$ (теорема~\ref {thm:2'}), ана\-ло\-гич\-ную
полному инварианту Ю.М.~Бурмана~\cite{B,B2} и доказываем равенство $|\pi_0(F)|=\infty$.

Близкая задача была решена С.В.~Матвеевым~\cite{Kmsb} (1997),
Х.~Цишангом~\cite {Z} (1998), В.В.\ Шарко~\cite{SH} (1998) и С.И.\
Максименко~\cite{Max2005} (2005). Матвеев и Цишанг доказали разными
методами линейную связность про\-стран\-ства $\tilde F=\tilde
F_{p,q,r}\supset F$ функций Морса на $M$, имеющих фиксированные
множества кри\-ти\-чес\-ких точек локальных минимумов и максимумов.
Другой близкий результат был получен Бурманом~\cite{B,B2}. Он изучал
пространство $F'$ гладких функций без критических точек на
некомпактной поверхности $M'$, локально постоянных на крае и имеющих
заданное поведение вблизи края. Для любой функции $f\in F'$ он
построил отображение $B_f\col F'\to H^1(M',\d M')$ (полный
топологический инвариант на пространстве $F'$) и доказал, что
индуцированное отображение $(B_f)_\#\col\pi_0(F')\to H^1(M',\d M')$
биективно. Гомотопический тип пространств функций с умеренными
особен\-нос\-тями на окружности изучался В.И.\ Арнольдом~\cite{A}.
Функции Морса на поверхностях изучались А.Т.\ Фоменко и Цишангом
\cite{FZ}, А.В.\ Бол\-си\-но\-вым и Фоменко~\cite{BF,BF1}, и
автором~\cite{K} в связи с задачей классификации
не\-вы\-рожд\-ен\-ных интегрируемых гамиль\-то\-но\-вых систем.
Различные вопросы классификации и топологии про\-странств функций
Морса на поверхностях исследовались также в работах~\cite
{Kul,Ma,KP2,KP1}.

Опишем основные результаты настоящей работы.

\begin{Not} \label{not:groups}
Пусть $\G=\G(M,\{x_1,\dots,x_{p+q+r}\})$ -- группа сохраняющих
ориентацию диффе\-о\-мор\-физ\-мов $M$, оставляющих неподвижными все
критические точки, пусть $\G^0$ -- компонента связности $\id_M$ в
$\G$, и $\G_f\subset\G$ -- множество диффеоморфизмов, сохраняющих
компоненту связности $F_f$ функции $f\in F$ в $F$ (в
$C^\infty$-топологиях на $\G$ и $F$, см.~\cite[\S4]{KP1}). Ниже
(определение~\ref {def:adm}) вводятся группа $\D_f^\abs$ {\it
абсолютно допустимых} и группа $\D_f$ {\it допустимых}
диффеоморфизмов для функции $f$ (отличные от понятия $f$-допустимого
диффеоморфизма $h\in\G_f$ из~\cite[\S6]{Max2005}). По теореме~\ref
{thm:2''} ниже, они являются нормальными подгруппами группы $\G_f$.
Так как группа $\G/\G^0$ дискретна, то подгруппы
$\D_f^\abs/\G^0\subset\D_f/\G^0\subset\G_f/\G^0$ и факторгруппы
$\G/\<\G^0\>$, $\G_f/\D_f$ и $\D_f/\D_f^\abs\subset\G_f/\D_f^\abs$
дискретны.
\end{Not}

Возникают следующие задачи:

1) Для заданного диффеоморфизма $h\in\G$ определить, принадлежат ли
функции $f$ и $fh$ одной компоненте связности $F_f$ пространства $F$
(т.е.\ принадлежит ли $h$ подгруппе $\G_f$). В частности, описать
пространство смежных классов $\G/\G_f\approx\pi_0(F)$ и определить,
является ли оно конечным.

2) Для заданного диффеоморфизма $h\in\G$ или $\G_f$ определить,
является ли он допустимым (абсолютно допустимым) для функции $f$
(т.е.\ принадлежит ли подгруппам $\D_f$ и $\D_f^\abs$). В частности,
подтвердить или опровергнуть гипотезу М.~Басмановой о совпадении
подгрупп $\D_f^\abs\subset\D_f\subset\G_f$.

3) Описать конечные множества порождающих элементов факторгрупп
$\G_f/\G^0$ и $\G_f/\D_f$.

В данной работе с помощью {\it числа вращения}, введенного
Б.~Рейнхартом~\cite {R}, получены частичные решения первых двух
задач, а с помощью {\it комплексов функций Морса} -- решение третьей
задачи:

1) Построена сюръекция $\pi_0(F)\approx\G/\G_f\to\ZZ^{p+r-1}$
(теорема~\ref {thm:2'}). В частности, $|\pi_0(F)|=\infty$.

2) Построена сюръекция $\G/\D_f^\abs\to\ZZ_2^{q-1}$, которая
индуцирует эпиморфизм $\G_f/\D_f^\abs\to\ZZ_2^{q-1}$, а при $M\ne
S^2$ -- эпиморфизм $\D_f/\D_f^\abs\to\ZZ_2^{q-1}$ (теорема~\ref
{thm:2''}). В частности, при $q\ge2$ мы получаем опровержение
$\D_f^\abs\subsetneq\G_f$ гипотезы Басмановой.

3) Определен конечный связный полиэдральный комплекс $K=K_{p,q,r}$,
ассо\-ци\-и\-ро\-ван\-ный с про\-стран\-ст\-вом $F$ (теорема~\ref
{thm:1}). Построены эпиморфизм $\mu\col\pi_1(K)\to\G_f/\D_f$ и
конечные множества порождающих элементов групп $\G_f/\G^0$ и
$\G_f/\D_f$ в терминах 2-остова комплекса $K$ (теоремы~\ref
{thm:K},~\ref {thm:BK}).

В статье также исследовано, какие из групп цепочки
$\G^0\subset\D_f^\abs\subset\D_f\subset\G_f\subset\G$ совпадают (см.\
следствие), кроме случая $\D_f\subset\G_f$ при $M\ne S^2$, $q\ge4$.
При $M=S^2$ доказаны оценки
$q-1\le\rank(\G_f/\D_f)\le\rank(\pi_1(K))$, $\rank(\G/\<\G_f\>)\ge
p+r-1$, а при $M\ne S^2$ оценки $\rank(\G_f/\D_f)\le\rank(\pi_1(K))$,
$\rank(\D_f/\D_f^\abs)\ge q-1$ (следствие и теорема~\ref {thm:K}),
где ранг группы есть минимальное количество порождающих элементов.
Отсюда $\D_f=\G_f$, если $\pi_1(K)=1$. Поэтому $\D_f=\G_f$ в случае
$M\ne S^2$, $q\le3$ (так как комплексы $K=K_{1,2,1}\approx[0,1]$ и
$K=K_{1,3,2}\sim\bigvee\limits_{7}S^2$ односвязны).

\medskip
{\bf 2. Топологический инвариант на пространстве $F$,
$\G_f$-инвариант на пространстве $\G$.} Обозначим через $\H$
подгруппу в $\G$, порожденную $\G^0$ и скручиваниями Дэна~\cite
{Dehn} вокруг разбивающих кривых (``ядро Джонсона''~\cite{J}). Она
является нормальной.

\begin{Thm} [$\G_f$-инвариант $B_f$ на пространстве $\G$] \label{thm:2'}
Пусть $q\ge1$ и $f\in F$. Имеется сюръекция
$B_f\col\G\to\ZZ^{p+r-1}$, ограничение которой на любой смежный класс
$\G_fh$, $h\in\G$, постоянно. Огра\-ни\-че\-ние
$B_f|_\H\col\H\to\ZZ^{p+r-1}$ не зависит от функции $f$ и является
эпиморфизмом.
Скручивание Дэна вокруг границы любого диска, содержащего ровно
$k\ge0$ седловых критических точек и $\ell\not\in\{0,k+1,p+r\}$
критических точек локальных минимумов и максимумов,
не принадлежит подгруппе $\G_f\cap\H\subset\ker(B_f|_\H)$ (т.е.\ не
сохраняет компоненту $F_f$ функции $f$ в $F$). В частности,
$|\pi_0(F)|=\infty$, $\G_f\subsetneq\G$ и имеется сюръекция
$\pi_0(F)\approx\G/\G_f\to\ZZ^{p+r-1}$. Если $M=S^2$, то $\H=\G$, и
$B_f$ определяет эпиморфизм $\G/\<\G_f\>\to\ZZ^{p+r-1}$, не зависящий
от $f$.
\end{Thm}

\medskip
{\bf 3. Допустимые диффеоморфизмы и $\D_f^\abs$-инвариант на
пространстве $\G$.}
\begin{Def} \label {def:adm}
Диффеоморфизм $h\in\G$ назовем {\it допустимым для функции} $f\in F$,
если имеются такие функции $f_1,\dots,f_N\in F_f$ и диффеоморфизмы
$h_1,\dots,h_N\in\G$, что $f_i=f_ih_i$ и $h\in h_1\dots h_N\G^0$.
Если каждый $h_i$ -- скручивание Дэна вокруг связной компоненты
кривой $f_i^{-1}(a_i)$, где $a_i$ -- не\-кри\-ти\-чес\-кое значение
функции $f_i$, то диффеоморфизм $h$ назовем {\it абсолютно допустимым
для} $f$. Абсолютно допустимые и допустимые диффеоморфизмы для
функции $f\in F$ образуют подгруппы $\D_f^\abs$ и $\D_f$ группы $\G$
(см.\ обозначение~\ref {not:groups}). Ясно, что $\G^0\subset
\D_f^{\abs} \subset \D_f \subset \G_f \subset \G$.
 \end{Def}

\begin {Exs} \label {examples:admiss}
(А) Простая замкнутая кривая на $M$ называется {\it допустимой} \cite
[\S6]{Max2005} для функции Морса $f\in F$, если она является
компонентой связности линии уровня $g^{-1}(a)$ некоторой функции
$g\in F_f$. Скручивание Дэна вокруг такой кривой -- это абсолютно
допустимый диффеоморфизм для $f$.

(Б) Другой пример допустимого диффеоморфизма показан на рис.~\ref
{fig:twist'}. Как и в примере (А), этот диффеоморфизм $h=h_{ij}$ сохраняет функцию $g\in F_f$,
однако он совпадает с тождественным в окрестностях всех критических
точек $x_1,\dots,x_{p+q+r}$ кроме двух седловых точек $x_i$ и $x_j$,
в которых $dh(x_i)=-\id$ и $dh(x_j)=-\id$. Такой диффеоморфизм
существует для любой поверхности $M\ne S^2$, а при $M=S^2$ -- нет. Он
не является абсолютно допустимым для $f$, согласно теореме~\ref
{thm:2''}(Б) ниже.

 \begin{figure}[htbp]
\setlength{\unitlength}{9pt}
\begin{center}
\begin{picture}(35,11)(-2,2)
\put(-2,0){
\qbezier[200](2,3.7)(2.9,2.7)(2.3,1.8)  
\qbezier[200](2,3.7)(-.6,7)(2,10.3)    
\qbezier[200](2,10.3)(2.9,11.3)(2.3,12.2)  
\qbezier[200](8,10.3)(10.6,7)(8,3.7)    
\qbezier[200](8,3.7)(7.1,2.7)(7.7,1.8)  
\qbezier[200](8,10.3)(7.1,11.3)(7.7,12.2)  
\qbezier[20](2.6,2.2)(5,3.25)(7.4,2.2) 
\qbezier[150](2.6,2.2)(2,1.7)(2.6,1.2)  
\qbezier[150](7.4,2.2)(8,1.7)(7.4,1.2)  
\qbezier[150](2.6,1.2)(5,0.15)(7.4,1.2) 
\put(0,10.6){
\qbezier[150](2.6,2.2)(5,3.25)(7.4,2.2) 
\qbezier[150](2.6,2.2)(2,1.7)(2.6,1.2)  
\qbezier[150](7.4,2.2)(8,1.7)(7.4,1.2)  
\qbezier[150](2.6,1.2)(5,0.15)(7.4,1.2) 
}
\thicklines
\qbezier[300](3,1)(3.8,3)(2.9,4.5)
\qbezier[300](2.9,4.5)(1.5,7)(2.8,9.3) 
\qbezier[300](2.8,9.3)(3.3,10.2)(3,11.6)
\qbezier[300](7,1)(6.2,3)(7.1,4.5)
\qbezier[300](7.1,4.5)(8.5,7)(7.2,9.3) 
\qbezier[300](7.2,9.3)(6.7,10.2)(7,11.6)
\thinlines
\put(5,7){
\qbezier[50](.3,-1)(-.8,0)(.3,1)
\qbezier[50](.1,-.8)(.8,0)(.1,.8)
\put(.1,-.8){\put(-.43,-.13){\tiny$\bullet$} \put(-1.5,-.6){$x_i$}}
\put(.1,.8){\put(-.43,-.23){\tiny$\bullet$} \put(-1.5,.5){$x_j$}}
}
}
\put(9,6){$\stackrel{\mbox{$h_{ij}$}}{\mbox{\LARGE$\longrightarrow$}}$}
\put(13.5,0){
\qbezier[200](2,3.7)(2.9,2.7)(2.3,1.8)  
\qbezier[200](2,3.7)(-.6,7)(2,10.3)    
\qbezier[200](2,10.3)(2.9,11.3)(2.3,12.2)  
\qbezier[200](8,10.3)(10.6,7)(8,3.7)    
\qbezier[200](8,3.7)(7.1,2.7)(7.7,1.8)  
\qbezier[200](8,10.3)(7.1,11.3)(7.7,12.2)  
\qbezier[20](2.6,2.2)(5,3.25)(7.4,2.2) 
\qbezier[150](2.6,2.2)(2,1.7)(2.6,1.2)  
\qbezier[150](7.4,2.2)(8,1.7)(7.4,1.2)  
\qbezier[150](2.6,1.2)(5,0.15)(7.4,1.2) 
\put(8,1.5) {\scriptsize$g^{-1}(a)$}
\put(0,10.6){
\qbezier[150](2.6,2.2)(5,3.25)(7.4,2.2) 
\qbezier[150](2.6,2.2)(2,1.7)(2.6,1.2)  
\qbezier[150](7.4,2.2)(8,1.7)(7.4,1.2)  
\qbezier[150](2.6,1.2)(5,0.15)(7.4,1.2) 
\put(8,1.5) {\scriptsize$g^{-1}(b)$}
}
\thicklines
\qbezier[300](3,1)(3,3.5)(5,3.5)
\qbezier[300](5,3.5)(6.6,3.5)(7.1,4.5)
\qbezier[300](7.1,4.5)(8.2,6.4)(7.6,8.3) 
\qbezier[300](7.6,8.3)(7.1,9.4)(4,9.1)
\qbezier[300](4,9.1)(2.5,9.1)(3,11.6)
\qbezier[300](7,1)(6.5,3.4)(7.8,3.8)
\qbezier[80](7.8,3.8)(8.1,3.9)(8.37,4.2)
\qbezier[6](8.37,4.2)(8.7,4.7)(7,5)
\qbezier[13](7,5)(5,5.3)(3,5)
\qbezier[6](1.63,4.2)(1.3,4.7)(3,5)
\qbezier[80](2.2,3.8)(1.9,3.9)(1.63,4.2)
\qbezier[100](2.2,3.8)(3.5,3.6)(2.9,4.5)
\qbezier[300](2.9,4.5)(1.8,6.4)(2.4,8.3) 
\qbezier[100](2.4,8.3)(2.6,8.9)(2.2,9.1)
\put(0,5.3){
\qbezier[80](8.2,4)(8.1,3.9)(8.37,4.2)
\qbezier[6](8.37,4.2)(8.7,4.7)(7,5)
\qbezier[13](7,5)(5,5.3)(3,5)
\qbezier[6](1.63,4.2)(1.3,4.7)(3,5)
\qbezier[80](2.2,3.8)(1.9,3.9)(1.63,4.2)
}
\qbezier[100](8.2,9.3)(7,9)(7,11.6)
\thinlines
\put(5,7){
\qbezier[50](.3,-1)(-.8,0)(.3,1)
\qbezier[50](.1,-.8)(.8,0)(.1,.8)
\put(.1,-.8){\put(-.43,-.13){\tiny$\bullet$} \put(-1.5,-.6){$x_i$}}
\put(.1,.8){\put(-.43,-.23){\tiny$\bullet$} \put(-1.5,.5){$x_j$}}
}
\put(12.5,6.3){$\stackrel{g}{\longrightarrow}$}
\put(18,1){\vector(0,1){13}}
\put(18,1.7){\put(-.43,-.18){\tiny$\bullet$} \put(.5,-.2){$a$}}
\put(18,12.3){\put(-.43,-.18){\tiny$\bullet$} \put(.5,-.2){$b$}}
\put(18,6.2){\put(-.43,-.18){\tiny$\bullet$} \put(.5,-.2){$c_i$}}
\put(18,7.8){\put(-.43,-.18){\tiny$\bullet$} \put(.5,-.2){$c_j$}}
}
\end{picture}
\end{center}
\setlength{\unitlength}{1pt}
 \fcaption {Допустимый, но не абсолютно допустимый диффеоморфизм $h_{ij}$.}
 \label {fig:twist'}
 \end {figure}

Группа $\D_f$ допустимых диффеоморфизмов порождена диффеоморфизмами
из примеров~\ref {examples:admiss}(А,Б).
 \end {Exs}

\begin{Thm} [$\D_f^\abs$-инвариант $B_f^\abs$ на пространстве $\G$] \label{thm:2''}
Пусть $q\ge1$ и $f\in F$. Имеется сюръекция
$B_f^\abs\col\G\to\ZZ_2^{q-1}$, ограничение которой на любой смежный
класс $\D_f^\abs h$, $h\in\G$, постоянно. Подгруппы $\D_f^\abs$ и
$\D_f$ являются нормальными в $\G_f$, и выполнены следующие условия:

{\rm (A)} Ограничение $B_f^\abs|_H\col H\to\ZZ_2^{q-1}$ на любую из
трех подгрупп $H\in\{\G_f,\H,\G_f\cap\H\}$ является эпиморфизмом, при
$H=\H$ не зависит от $f$, и $\D_f^\abs\subset\ker(B_f^\abs|_{\G_f})$.
При $q\ge 2$ для любой пары седловых критических точек
скручивание Дэна вокруг границы некоторого диска (зависящего от $f$),
содержащего эти две точки и не содержащего других критических точек,
принадлежит $\G_f\setminus\D_f^\abs$ (т.е.\ сохраняет компоненту
$F_f$ функции $f$ в $F$, но не является абсолютно допустимым
диффе\-о\-мор\-физ\-мом для $f$). В частности, ограничение
$B_f^\abs|_{\G_f}$ индуцирует эпиморфизм
$\G_f/\D_f^\abs\to\ZZ_2^{q-1}$, и поэтому $\D_f^\abs\subsetneq\G_f$
при $q\ge2$. Если $M=S^2$ и $q\ge2$, то
$\D_f^\abs=\D_f\subsetneq\G_f$.

{\rm (Б)} Если $M\ne S^2$, то ограничение
$B_f^\abs|_{\D_f}\col\D_f\to\ZZ_2^{q-1}$ является эпиморфизмом,
инду\-ци\-ру\-ющим эпиморфизм $\D_f/\D_f^\abs\to\ZZ_2^{q-1}$, причем
$\D_f^\abs\subsetneq\D_f\subset\G_f$ и допустимый для функции $f$
диффе\-о\-мор\-физм, показанный на рис.~$\ref {fig:twist'}$, не
является абсолютно допустимым для $f$.
\end{Thm}

Если у функции $f_1\in F_f$ ровно $q\ge1$ седловых критических
значений, то на $M$ имеются $q+g-1$ окружностей, являющихся
компонентами линий уровня функции $f_1$, и таких что подгруппа группы
$\D_f^\abs/\G^0$, порожденная скручиваниями Дэна вокруг этих
окружностей, изоморфна $\ZZ^{q+g-1}$.

\begin{Cor}
{\rm (A)} Пусть $M=S^2$. Если количество седел $q\ge2$, то имеется
цепочка четырех групп
$\G^0\subsetneq\D_f^\abs=\D_f\subsetneq\G_f\subsetneq\G$, в которой
все множества смежных классов бесконечны и допускают мономорфизм
$\ZZ^{q-1}\rightarrowtail\D_f^\abs/\G^0$ и эпиморфизмы
$\G_f/\D_f\to\ZZ_2^{q-1}$, $\G/\<\G_f\>\to\ZZ^{p+r-1}$. Если $q=1$,
то имеются две группы $\G^0=\D_f^\abs=\D_f=\G_f\subsetneq\G$ с
бесконечной факторгруппой
$\G/\G^0\cong\pi_1(S^2\setminus\{x_2,x_3,x_4\},x_1)\cong F_2$, где
$F_2$ -- свободная группа ранга $2$.

{\rm (Б)} Если $M\ne S^2$, то имеется цепочка пяти групп
$\G^0\subsetneq\D_f^\abs\subsetneq\D_f\subset\G_f\subsetneq\G$, в
которой все множества смежных классов (за исключением, быть может,
$\G_f/\D_f$) бесконечны и допускают мономорфизм
$\ZZ^{q+g-1}\rightarrowtail\D_f^\abs/\G^0$, эпиморфизм
$\D_f/\D_f^\abs\to\ZZ_2^{q-1}$ и сюръекцию $\G/\G_f\to\ZZ^{p+r-1}$.
\end{Cor}

\PProof{теорем~\ref{thm:2'} и~\ref{thm:2''}} Шаг 1. В данном
доказательстве под {\it кривой} понимается гладкое компактное (не
обязательно связное) ориентированное 1-мерное подмногообразие
$\alpha\subset M$, край которого есть пересечение множества $\alpha$
с множеством критических точек $x_1,\dots,x_{p+q+r}$. Пусть
 \begin{equation} \label {eq:gamma:i}
\gamma_i\colon\ [0,1]\to M, \quad 1\le i\le p+q+r-1,
 \end{equation}
-- кривая из точки $\gamma_i(0)=x_{p+q+r}$ в точку $\gamma_i(1)=x_i$.
Фиксируем на $M$ риманову метрику.

\begin{Def}
Для любой такой кривой $\gamma\col [0,1]\to M$ и любой функции $f\in
F$ обозначим через $w_f(\gamma)$ вещественное число, равное ``полному
количеству оборотов'' касательного вектора ${d\gamma\over dt}(t)$
вокруг нуля по отношению к ортогональному реперу в $T_{\gamma(t)}M$,
содержащему вектор $\grad f(\gamma(t))$, $0<t<1$. Для несвязной
кривой $\gamma\subset M$ определим $w_f(\gamma)$ равным сумме чисел,
отвечающих ее компонентам. Назовем $w_f(\gamma)$ {\it числом вращения
кривой $\gamma$ по отношению к функции} $f$. Оно совпадает с {\it
числом вращения кривой $\gamma$ по отношению к векторному полю}
$\grad f$ (см.~\cite[\S2]{R}, \cite[определения (1.1)]{Ch} или
\cite[\S3.2]{B2}). Для замкнутой кривой $\gamma$ число $w_f(\gamma)$
целое и не меняется при деформациях функции $f\in F$
(см.~\cite[\S2]{R}, \cite[леммы (5.1) и (5.2)]{Ch}, \cite[\S3.1,
утверждение 5]{B2}).
\end{Def}

Аналогично~\cite[\S2]{R}, определим {\it различающее число} кривой
$\gamma$ по отношению к функциям $f,fh$:
 $$
\d_hw_f(\gamma):=w_f(h\gamma)-w_f(\gamma)=w_{fh}(\gamma)-w_f(\gamma)
= (w_{fh}-w_f)(\gamma), \quad h\in\G.
 $$

Отметим некоторые свойства чисел $w_f(\gamma)$ и $\d_hw_f(\gamma)$.
Для любой пары $h_1,h_2\in\G$ выполнено
 \begin{equation} \label {eq:*}
 \d_{h_1h_2}w_f(\gamma) = \d_{h_1}w_f(\gamma) + \d_{h_2}w_{fh_1}(\gamma),
 \end{equation}
поскольку
$\d_{h_1h_2}w_f(\gamma)=(w_{fh_1h_2}-w_f)(\gamma)=(w_{fh_1h_2}-w_{fh_1}+w_{fh_1}-w_f)(\gamma)
 =(\d_{h_2}w_{fh_1}+\d_{h_1}w_f)(\gamma)$.
Если $s_i$ -- маленькая окружность вокруг точки $x_i$,
ориентированная ``против часовой стрелки'', то
 \begin{equation} \label {eq:A0}
w_f(s_i)=1-\ind_{x_i}(\grad f), \quad 1\le i\le p+q+r.
 \end{equation}
Таким образом, $w_f(s_i)$ всегда четно, так как $w_f(s_i)=0$ для
точек $x_i$ локальных минимумов и максимумов ($q<i\le p+q+r$) и
$w_f(s_i)=2$ для седловых точек $x_i$ ($1\le i\le q$). Более общо,
для любой (не обязательно связной) разбивающей кривой $\alpha=\d N$,
где $N\subset M$, выполнено
 \begin{equation} \label {eq:DeltaAN}
w_f(\d N) = \chi(N) - \sum_{x_i\in N} \ind_{x_i}(\grad f),
 \end{equation}
где кривая $\d N$ ориентирована так, что $N$ ``находится слева'' (это
выводится из~(\ref {eq:A0}) приклеиванием дисков к компонентам $\d N$
и продолжением векторного поля $\grad f$ внутрь каждого диска с одной
особой точкой, см.~\cite[лемма (5.7)]{Ch}). Для любой кривой $\gamma$
и любой связной замкнутой кривой $\alpha$
 \begin{equation} \label {eq:DeltaA0}
\d_{t_\alpha^k}w_f(\gamma) = k\la\alpha,\gamma\ra w_f(\alpha), \qquad
k\in\ZZ,
 \end{equation}
где $\la\alpha,\gamma\ra$ -- индекс пересечения кривых $\alpha$ и
$\gamma$, $t_\alpha$ -- скручивание Дэна вокруг $\alpha$.

Согласно (\ref {eq:A0}), (\ref {eq:DeltaA0}) и построению (\ref
{eq:gamma:i}) кривых $\gamma_i$, для любого $j\in[1,p+q+r-1]$
выполнено
 \begin{equation} \label {eq:DeltaA}
 \d_{t_{s_j}}w_f(\gamma_i) = 2\delta_{ij}, \quad 1\le i\le q;
 \qquad
 \d_{t_{s_j}}w_f(\gamma_i) = 0, \quad q < i\le p+q+r-1.
 \end{equation}
Для любого $h\in\G$ выберем диффеоморфизм $\tilde h\in h\G^0$,
ограничение которого на малую окрестность $U$ множества точек
$x_1,\dots,x_{p+q+r}$ совпадает с $\id_M$. Ясно, что число
$\d_{\tilde h}w_f(\gamma_i)$ целое и сохраняется при деформациях
функции $f$ в $F$. При этом, в силу~(\ref {eq:*}) и~(\ref
{eq:DeltaA}), значения $\d_{\tilde h}w_f(\gamma)\mod2\in\ZZ_2$, $1\le
i\le q$, и $\d_{\tilde h}w_f(\gamma)\in\ZZ$, $q<i\le p+q+r-1$, не
зависят от выбора диффеоморфизма $\tilde h$. Для любых функции $f\in
F$ и набора кривых (\ref {eq:gamma:i}) определим отображения $B_f$ и
$B_f^\abs$ формулами
 $$
 B_f^\abs\col \G\to \ZZ_2^{q-1}, \qquad B_f^\abs(h)
 := (\d_{\tilde h}w_f(\gamma_1)\mod2,\dots, \d_{\tilde h}w_f(\gamma_{q-1})\mod2);
 $$
 $$
 B_f\col \G\to \ZZ^{p+r-1}, \qquad B_f(h)
 := (\d_{\tilde h}w_f(\gamma_{q+1}),\dots, \d_{\tilde h}w_f(\gamma_{p+q+r-1})).
 $$

В силу (\ref {eq:*}), для любых $h_1,h_2\in\G$ выполнены равенства
 \begin{equation} \label {eq:**}
 B_f(h_1h_2) = B_f(h_1) + B_{fh_1}(h_2) , \quad
 B_f^\abs(h_1h_2) = B_f^\abs(h_1) + B_{fh_1}^\abs(h_2) .
 \end{equation}
Поэтому для любых $h_1\in\G_f$ и $h_2\in\G$ выполнены (в силу
$B_{fh_1}=B_f$ и $B_{fh_1}^\abs=B_f^\abs$) равенства
 \begin{equation} \label {eq:***}
 B_f(h_1h_2) = B_f(h_1) + B_f(h_2), \quad
 B_f^\abs(h_1h_2) = B_f^\abs(h_1) + B_f^\abs(h_2).
 \end{equation}

Шаг 2. Докажем равенство $B_f(\G_f h_2)=B_f(h_2)$ для любого
$h_2\in\G$. Сначала докажем равенство $B_f(\G_f)=0$. Для любого
$h\in\G_f$ рассмотрим число $\d_{\tilde h}w_f(\gamma_{i})=w_{f\tilde
h}(\gamma_i)-w_f(\gamma_i)$, $q<i<p+q+r$. Пусть $U'\subset U$ --
малая окрестность множества $\{x_{q+1},\dots,x_{p+q+r}\}$ точек
локальных минимумов и максимумов. Тогда любой путь $f_t$ в $F$ со
свойством $f_0|_U=f_1|_U$ гомотопен в классе путей с фиксированными
концами в пространстве $F$ такому пути $\tilde f_t$, что $\tilde
f_t|_{U'}=f_0|_{U'}$ при любом $t\in[0,1]$. Из $h\in\G_f$ имеем
$f\tilde h\in F_f$, поэтому существует путь $f_t$ в $F$, такой что
$f_0=f$, $f_1=f\tilde h$ и $f_t|_{U'}=f|_{U'}$ при любом $t\in[0,1]$.
Разность $w_{f_t}(\gamma_i)-w_f(\gamma_i)$ целая при любом $t$ (так
как концы кривой $\gamma_i$ содержатся в $U'$), а значит, постоянна и
равна $w_{f\tilde
h}(\gamma_i)-w_f(\gamma_i)=w_{f_0}(\gamma_i)-w_f(\gamma_i)=0$.
Поэтому $B_f(h)=0$ и $B_f(\G_f)=0$. С учетом (\ref{eq:***}), это дает
$B_f(\G_f h_2)=B_f(\G_f)+B_f(h_2)=B_f(h_2)$.

Докажем равенство $B_f^\abs(\D_f^\abs h_2)=B_f^\abs(h_2)$ для любого
диффеоморфизма $h_2\in\G$. Заметим, что $w_f(\alpha)=0$ для любой
допустимой кривой $\alpha$ для $f$ (см.\ определение~\ref {def:adm}).
В силу (\ref{eq:DeltaA0}), это дает равенство
$\d_{t_\alpha^k}w_f(\gamma_i)=0$ при $1\le i<p+q+r$, $k\in\ZZ$,
откуда $B_f^\abs(t_\alpha^k)=0$. С учетом (\ref {eq:***}), для любого
$h_2\in\G$ выполнено $B_f^\abs(t_\alpha^kh_2)=B_f^\abs(h_2)$, откуда
индукцией получаем $B_f^\abs(\D_f^\abs h_2)=B_f^\abs(h_2)$.

Шаг 3. Докажем, что отображения $B_f^\abs|_{\G_f}$, $B_f^\abs|_\H$ и
$B_f|_\H$ являются гомоморфизмами, причем второй и третий не зависят
от функции $f\in F$. Первое отображение является гомоморфизмом в силу
(\ref {eq:***}). В силу (\ref {eq:DeltaAN}), для любой связной
разбивающей кривой $\alpha=\d N$ число $w_f(\alpha)$ не зависит от
$f$. С учетом (\ref {eq:DeltaA0}), для любого $k\in\ZZ$ число
$\d_{t_\alpha^k}w_f(\gamma)=k\la\alpha,\gamma\ra w_f(\alpha)$ тоже не
зависит от $f$. Поэтому $B_f(t_\alpha^k)$ не зависит от $f$. Отсюда и
из (\ref {eq:**}) получаем, что $B_f(h_1t_\alpha^k) = B_f(h_1) +
B_{fh_1}(t_\alpha^k) = B_f(h_1) + B_f(t_\alpha^k)$ для любого
$h_1\in\G$. Поэтому $B_f|_\H$ -- гомоморфизм и не зависит от $f$;
аналогичное верно для $B_f^\abs|_\H$.

Шаг 4. Покажем, что гомоморфизмы $B_f^\abs|_{\G_f\cap\H}$ и $B_f|_\H$
являются эпиморфизмами. Это следует из следующего факта. Для любых
функции $f\in F$ и числа $i\ne q$, $1\le i<p+q+r$ (точнее, $i<q$ для
$B_f^\abs|_{\G_f\cap\H}$ и $i>q$ для $B_f|_\H$) можно построить
замкнутую кривую $s_{iq}=s_i\# s_q\subset M$, являющуюся ``связной
суммой'' маленьких окружностей $s_i$ и $s_q$ вокруг критических точек
$x_i$ и $x_q$ и такую, что скручивание Дэна $t_{s_{iq}}$ вокруг кривой $s_{iq}$
обладает следующими свойствами:
 \begin{enumerate}
 \item
$t_{s_{iq}}\in\H$, а в случае $1\le i<q$ выполнено
$t_{s_{iq}}\in\G_f$ (т.е.\ функция $ft_{s_{iq}}$ принадлежит
компоненте связности $F_f$ функции $f$ в пространстве $F$);
 \item
в случае $1\le i<q$ элемент $B_f^\abs(t_{s_{iq}})$ совпадает с $i$-ым
элементом канонического базиса группы $\ZZ_2^{q-1}$, а в случае
$q<i<p+q+r$ элемент $B_f(t_{s_{iq}})$ совпадает с $(i-q)$-ым
элементом канонического базиса группы $\ZZ^{p+r-1}$ (поэтому
$B_f^\abs(\G_f\cap\H)=\ZZ_2^{q-1}$ и $B_f(\H)=\ZZ^{p+r-1}$).
 \end{enumerate}
Первая часть пункта 1 следует из определения группы $\H$ (так как
$s_{iq}$ -- связная разбивающая кривая). Пункт 2 следует из~(\ref
{eq:DeltaA0}) и~(\ref {eq:DeltaAN}), так как (для любого $j\ne q$,
$1\le j<p+q+r$) $\d_{t_{s_{iq}}}w_f(\gamma_j)=\la s_{iq},\gamma_j\ra
w_f(s_{iq})$ равно $\la s_{iq},\gamma_j\ra\cdot3=3\delta_{ij}$ при
$1\le i<q$ и равно $\la s_{iq},\gamma_j\ra\cdot1=\delta_{ij}$ при
$q<i<p+q+r$. Заменяя окружность $s_{iq}$ на границу диска $D\subset
M$, содержащего $k$ седловых и $\ell\not\in\{0,k+1,p+r\}$ минимаксных
критических точек, а кривую $\gamma_j$ на любую кривую $\gamma$,
ведущую из точки минимакса снаружи $D$ в точку минимакса в $D$,
из~(\ref {eq:DeltaA0}) и~(\ref {eq:DeltaAN}) аналогично получаем
$\d_{t_{\d D}}w_f(\gamma)=\la \d D,\gamma\ra w_f(\d
D)=1\cdot(1+k-\ell)\ne0$, откуда $t_{\d D}\not\in\G_f$ (так как
$\d_{\tilde h}w_f(\gamma)=0$ для любого $h\in\G_f$, см.\ шаг~2).

 \begin{figure}[htbp]
\setlength{\unitlength}{12pt}
\begin{center}
\begin{picture}(41,15)(0,0)
\put(4,12){
\thicklines \put(-4,3){\line(1,0){2.5}} \put(4,3){\line(-1,0){2.5}}
\put(-1.5,-3){\line(1,0){3}}
\thinlines \qbezier[200](-2.75,3)(-2.2,1.6)(-1.2,0)
\put(-1.2,0){\put(-.386,-.223){$\bullet$} \put(-1.4,-.2){$x_i$}}
\qbezier[200](-1.2,0)(-.5,-1.5)(-.5,-3)
\put(-.5,-2.9){\vector(0,-1){.1} \put(-1.1,1.2){$\al$}}
\qbezier[200](2.75,3)(2.2,1.6)(1.2,0)
\put(1.2,0){\put(-.386,-.223){$\bullet$} \put(.3,-.2){$x_q$}}
\qbezier[200](1.2,0)(.5,-1.5)(.5,-3) \put(.5,-2.9){\vector(0,-1){.1}
\put(.3,1.2){$\be$}} \put(4.5,-.5){\huge $\rightsquigarrow$}
\put(-.3,-3.2){\rotatebox{-90}{\LARGE $\longrightarrow$}}
\put(.3,-4.6){$t_{s_{iq}}$} }
\put(15,12){
\thicklines \put(-4,3){\line(1,0){2.5}} \put(4,3){\line(-1,0){2.5}}
\put(-1.5,-3){\line(1,0){3}}
\thinlines
\qbezier[150](2.3,3)(1,.7)(-1.2,0)    
\put(-1.2,0){\put(-.386,-.223){$\bullet$}}
\qbezier[200](-2.75,3)(-2.2,1.6)(-1.2,0)
\qbezier[200](2.75,3)(2.2,1.6)(1.2,0) 
\put(1.2,0){\put(-.386,-.223){$\bullet$}}
\qbezier[200](1.2,0)(.5,-1.5)(.5,-3)
\put(4.5,-.5){\huge $\rightsquigarrow$} }
\put(26,12){
\thicklines \put(-4,3){\line(1,0){2.5}} \put(4,3){\line(-1,0){2.5}}
\put(-1.5,-3){\line(1,0){3}}
\thinlines
\qbezier[150](2.3,3)(1,.7)(-1.2,0)    
\put(-1.2,0){\put(-.386,-.223){$\bullet$}}
\qbezier[200](-2.5,3)(-2.2,1.6)(-1.2,0)
\qbezier[200](1.2,0)(.5,-1.5)(.5,-3)
\put(1.2,0){\put(-.386,-.223){$\bullet$}}
\qbezier[150](1.2,0)(0,-1.2)(-1.3,-.9)        
\qbezier[50](-1.3,-.9)(-2,-.6)(-2.3,.3)    
\qbezier[150](-2.3,.3)(-2.8,1.6)(-3,3)      
\put(4.5,-.5){\huge $\rightsquigarrow$} }
\put(37,12){
\thicklines \put(-4,3){\line(1,0){2.5}} \put(4,3){\line(-1,0){2.5}}
\put(-1.5,-3){\line(1,0){3}}
\thinlines \qbezier[50](2.5,3)(2.1,2)(1.4,1.7)
\qbezier[150](1.4,1.7)(-.2,1.1)(-1.2,0)    
\put(-1.2,0){\put(-.386,-.223){$\bullet$}}
\qbezier[150](-1.2,0)(1.3,1.5)(2.1,.7)     
\qbezier[150](2.1,.7)(2.5,.1)(1.9,-.9)       
\qbezier[150](1.9,-.9)(1.5,-1.5)(1.1,-2)   
\qbezier[150](1.1,-2)(.85,-2.3)(.7,-3)
\qbezier[150](1.2,0)(.1,-1.9)(0,-3) 
\put(1.2,0){\put(-.386,-.223){$\bullet$}}
\qbezier[150](1.2,0)(0,-1.2)(-1.3,-.9)        
\qbezier[50](-1.3,-.9)(-2,-.6)(-2.3,.3)
\qbezier[150](-2.3,.3)(-2.8,1.6)(-3,3)
\put(-.3,-3.5){\rotatebox{-90}{\huge $\rightsquigarrow$}} 
}
\put(37,3){
\thicklines \put(-4,3){\line(1,0){2.5}} \put(4,3){\line(-1,0){2.5}}
\put(-1.5,-3){\line(1,0){3}}
\thinlines \qbezier[50](3.2,3)(2.5,2)(1.4,1.7)
\qbezier[150](1.4,1.7)(-.2,1.1)(-1.2,0)    
\put(-1.2,0){\put(-.386,-.223){$\bullet$}}
\qbezier[150](-1.2,0)(1.3,1.5)(2.1,.7)     
\qbezier[150](2.1,.7)(2.5,.1)(1.9,-.9)       
\qbezier[150](1.9,-.9)(1.5,-1.5)(1.1,-2)   
\qbezier[150](1.1,-2)(.85,-2.3)(.7,-3)
\qbezier[150](-3,3)(-3.3,-.8)(-2.1,-1.1)  
\qbezier[150](1.2,0)(-.4,-1.9)(-2.1,-1.1) 
\put(1.2,0){\put(-.386,-.223){$\bullet$}}
\qbezier[150](1.2,0)(-.3,-1.2)(-1.8,-.7)        
\qbezier[150](-1.8,-.7)(-2.7,-.1)(-1.7,.9)    
\qbezier[150](-1.7,.9)(-1.1,1.5)(.4,2)    
\qbezier[50](.4,2)(1.5,2.2)(2.2,3) }
\put(26,3){
\thicklines \put(-4,3){\line(1,0){2.5}} \put(4,3){\line(-1,0){2.5}}
\put(-1.5,-3){\line(1,0){3}}
\thinlines
\qbezier[150](-1.2,0)(.3,1.2)(1.8,.7)        
\qbezier[150](1.8,.7)(2.7,.1)(2,-.7)    
\qbezier[150](2,-.7)(1.8,-.9)(1.4,-1.1)    
\qbezier[150](1.4,-1.1)(-.5,-2.2)(-2.7,-1.2)   
\qbezier[150](-2.7,-1.2)(-3.9,-.1)(-3.1,1.3) 
\qbezier[150](-3.1,1.3)(-2.6,2)(-2.75,3)
\put(-1.2,0){\put(-.386,-.223){$\bullet$}}
\qbezier[150](-1.2,0)(0,1.8)(2.2,1) 
\qbezier[150](2.2,1)(3.4,.1)(2.5,-1.1)  
\qbezier[50](2.5,-1.1)(1.9,-1.65)(1.6,-1.8)
\qbezier[50](1.6,-1.8)(1,-2.1)(.7,-3)
\qbezier[150](1.2,0)(0,-1.8)(-2.2,-1) 
\qbezier[150](-2.2,-1)(-3.3,-.1)(-2.5,1.1)  
\qbezier[150](-2.5,1.1)(-1.8,2.2)(-2,3)      
\put(1.2,0){\put(-.386,-.223){$\bullet$}}
\qbezier[150](1.2,0)(-.3,-1.2)(-1.8,-.7)        
\qbezier[150](-1.8,-.7)(-2.7,-.1)(-1.7,.9)    
\qbezier[150](-1.7,.9)(-1.1,1.5)(.4,2)    
\qbezier[50](.4,2)(1.5,2.2)(2.2,3) \put(4.5,.4){\rotatebox{180}{\huge
$\rightsquigarrow$}} }
\put(15,3){
\thicklines \put(-4,3){\line(1,0){2.5}} \put(4,3){\line(-1,0){2.5}}
\put(-1.5,-3){\line(1,0){3}}
\thinlines
\qbezier[150](-1.2,0)(.3,1.2)(1.8,.7)        
\qbezier[150](1.8,.7)(2.7,.1)(2,-.7)    
\qbezier[150](2,-.7)(1.8,-.9)(1.4,-1.1)    
\qbezier[150](1.4,-1.1)(-.5,-2.2)(-2.7,-1.2)   
\qbezier[150](-2.7,-1.2)(-3.9,-.1)(-3.1,1.3) 
\qbezier[150](-3.1,1.3)(-2.6,2)(-2.75,3)
\put(-1.2,0){\put(-.386,-.223){$\bullet$}}
\qbezier[150](-1.2,0)(0,1.8)(2.2,1) 
\qbezier[150](2.2,1)(3.4,.1)(2.5,-1.1)  
\qbezier[50](2.5,-1.1)(1.9,-1.65)(1.6,-1.7)
\qbezier[150](1.6,-1.7)(0,-2.3)(-.3,-3)
\qbezier[150](1.2,0)(0,-1.8)(-2.2,-1) 
\qbezier[150](-2.2,-1)(-3.3,-.1)(-2.5,1.1)  
\qbezier[150](-2.5,1.1)(-1.4,2.1)(.4,2.2)      
\qbezier[50](.4,2.2)(1.5,2.2)(2.2,3)
\put(1.2,0){\put(-.386,-.223){$\bullet$}}
\qbezier[150](1.2,0)(-.3,-1.2)(-1.8,-.7)        
\qbezier[150](-1.8,-.7)(-2.7,-.1)(-1.7,.9)    
\qbezier[150](-1.7,.9)(.5,2.4)(2.7,1.2)    
\qbezier[150](2.7,1.2)(3.9,.1)(3,-1.3)    
\qbezier[150](3,-1.3)(2.7,-1.8)(1.5,-2.3)
\qbezier[150](1.5,-2.3)(.6,-2.6)(.5,-3)
\put(4.5,.4){\rotatebox{180}{\huge $\rightsquigarrow$}} }
\put(4,3){
\thicklines \put(-4,3){\line(1,0){2.5}} \put(4,3){\line(-1,0){2.5}}
\put(-1.5,-3){\line(1,0){3}}
\thinlines \put(-1.2,0){\put(-.386,-.223){$\bullet$}}
\qbezier[150](1.2,0)(1.25,.7)(1.8,.5)
\qbezier[150](1.8,.5)(2.6,.1)(2.1,-.7)     
\qbezier[150](1.2,0)(.6,-1.98)(-2.1,-.7)  
\qbezier[150](-1.2,0)(-1.1,-.6)(-1.8,-.5)
\qbezier[150](-1.8,-.5)(-2.6,-.1)(-2.1,.7)  
\qbezier[150](-1.2,0)(-1.4,.8)(-.6,1.1)
\qbezier[150](-.6,1.1)(.4,1.5)(2.1,.7)     
\qbezier[150](2.1,.7)(3,.1)(2.4,-.9)       
\qbezier[150](2.1,-.7)(0,-2.3)(-2.4,-.9) 
\qbezier[150](-2.1,-.7)(-3,-.1)(-2.4,.9)    
\qbezier[150](-2.1,.7)(0,2.3)(2.4,.9)    
\qbezier[150](2.4,.9)(3.5,.1)(2.7,-1.1)     
\qbezier[150](2.4,-.9)(0,-2.7)(-2.7,-1.1)   
\qbezier[150](-2.4,-.9)(-3.5,-.1)(-2.7,1.1)  
\qbezier[150](-2.4,.9)(0,2.7)(2.7,1.1)      
\qbezier[150](2.7,1.1)(3.9,.1)(3,-1.3)    
\qbezier[150](3,-1.3)(2.7,-1.8)(1.5,-2.2)
\qbezier[150](1.5,-2.2)(.45,-2.35)(.5,-3)
\qbezier[150](2.7,-1.1)(1.6,-2.1)(.1,-2.15)   
\qbezier[150](.1,-2.15)(-.7,-2.3)(-.5,-3)
\qbezier[150](-2.7,-1.1)(-3.9,-.1)(-3,1.3) 
\qbezier[150](-3,1.3)(-1.9,1.9)(-2.75,3)
\qbezier[150](-2.7,1.1)(-.6,2.8)(2,1.8)      
\qbezier[150](2,1.8)(2.5,1.7)(2.75,3)
\put(1.2,0){\put(-.386,-.223){$\bullet$}}
\put(4.5,.4){\rotatebox{180}{\huge $\rightsquigarrow$}} }
\end{picture}
\end{center}
\setlength{\unitlength}{1pt}
 \fcaption {Реализация действия на $f$ скручивания вокруг двух седел гомотопией в $F$.}
 \label {fig:saddle:twist}
 \end {figure}

Осталось доказать вторую часть пункта 1. Мы построим требуемую кривую
$s_{iq}$, $1\le i<q$. Без ограничения общности считаем, что седловые
значения $f(x_i),f(x_q)$ превосходят остальные седловые значения
$f(x_j)$, $1\le j\le q-1$, $j\ne i$, и существует точка $x_k$
локального максимума, в которую входят сепаратрисы $\al$ и $\be$ поля
$\grad f$, выходящие из точек $x_i$ и $x_q$ соответственно. Пусть $D$
-- маленький круг вокруг $x_k$. Рассмотрим кривую $\al\cdot\be^{-1}$
и заменим ее часть $(\al\cdot\be^{-1})\cap D$ дугой окружности $\d
D$, не пересекающей две другие сепаратрисы, выходящие из точек $x_i$
и $x_q$ (существование такой дуги не ограничивает общности).
Рассмотрим связную сумму $s_{iq}=s_i\# s_q$ окружностей $s_i$ и $s_q$
по отношению к части полученной кривой между точками пересечения с
окружностями $s_i$ и $s_q$. Покажем, что существует путь из функции
$f$ в функцию $ft_{s_{iq}}$ в пространстве $F$ функций Морса. Этот
путь схематически изображен на рис.~\ref {fig:saddle:twist}.
Теорема~\ref {thm:2'} доказана.
%

Шаг 5. Покажем, что подгруппы $\D_f$ и $\D_f^\abs$ нормальны в
$\G_f$. Если $h_1\in\G_f$ (т.е.\ $fh_1\in F_f$) и диффеоморфизм
$d\in\G$ сохраняет функцию $fh_1$ (т.е.\ $fh_1d=fh_1$), то для любого
$h\in\G_f$ выполнено $(fh_1h^{-1})(hdh^{-1})=fh_1h^{-1}$, т.е.\
диффеоморфизм $hdh^{-1}$ сохраняет функцию $fh_1h^{-1}\in F_f$. Так
как группа $\D_f$ порождена $\G^0$ и всеми такими $d$ (или всеми
такими $hdh^{-1}$), то $h\D_f h^{-1}=\D_f$. Аналогично доказывается
равенство $h\D_f^\abs h^{-1}=\D_f^\abs$ (для этого в качестве $d$
рассматриваются лишь скручивания Дэна). Так как $h\in\G_f$ любой, то
подгруппы $\D_f$ и $\D_f^\abs$ нормальны в $\G_f$.

Шаг 6. Пусть $M\ne S^2$. Покажем, что $B_f^\abs(\D_f)=\ZZ_2^{q-1}$.
Рассмотрим допустимый, но не абсолютно допустимый для $f$
диффеоморфизм $h_{iq}\in\D_f\setminus\D_f^\abs\subset\G_f$,
показанный на рис.~\ref {fig:twist'}, при $1\le i<q$. Легко
проверяется, что $B_f^\abs(h_{iq})$ является $i$-ым элементом
стандартного базиса группы $\ZZ_2^{q-1}$. Поскольку $i$ любое и
$B_f^\abs|_{\G_f}$ -- гомоморфизм, то $B_f^\abs(\D_f)=\ZZ_2^{q-1}$.
Теорема~\ref {thm:2''} доказана.
 \qed

\medskip
{\bf 4. Эквивалентность и послойная эквивалентность функций Морса.} В
следствии описано, какие из соседних групп цепочки
$\G^0\subset\D_f^\abs\subset\D_f\subset\G_f\subset\G$ совпадают,
кроме случая $\D_f\subset\G_f$ при $M\ne S^2$. Наша дальнейшая цель
-- описать конечные множества поро\-жда\-ющих элементов факторгрупп
$\G_f/\D_f$ и $\G_f/\G^0$ в геометрических терминах.

\begin {Def} \label {def:equiv}
Функции Морса $f,g\in F$ назовем {\it подобными}, если они определяют
одно и то же разбиение поверхности $M$ на связные компоненты линий
уровня $f^{-1}(a)$ и $g^{-1}(b)$, а также один и тот же частичный
порядок на множестве седловых критических точек $x_1,\dots,x_q$
согласно значениям функции в этих точках; обозначим это следующим
образом: $f\approx g$. Если $f\approx gh$ для некоторого
диффеоморфизма $h\in\G$ (соответственно $h\in\G^0$), то функции $f,g$
назовем {\it эквивалентными} (соответственно {\it изотопными});
обозначим это через $f\sim g$ (соответственно $f\sim_\isot g$).
Классы эквивалентности и изотопности функции $f$ обозначим через
$[f]$ и $[f]_\isot$ соответственно.
 \end {Def}

Если в определении~\ref {def:equiv} не налагать условие о частичном
порядке на множестве седловых точек, получатся определения {\it
послойной подобности, послойной экви\-ва\-лент\-ности} и {\it
послойной изотоп\-ности}. Фоменко и Болсинов ввели комбинаторные
понятия {\it атома} и {\it молекулы} и доказали, что классы {\it
послойной эквивалентности} функций Морса на замкнутой поверхности
находятся во взаимно одно\-знач\-ном соответствии с молекулами таких
функций~\cite [Гл.~2, \S\S3--8, теорема~8]{BF1}. Ана\-ло\-гич\-но
вводятся понятия {\it нумерованного атома}, {\it нумерованной
молекулы} (с помощью нумерации вершин атомов согласно нумерации
критических точек $x_1,\dots,x_{p+q+r}$) и {\it оснащенной молекулы}
(с помощью частичного порядка из определения~\ref {def:equiv}) и
доказывается следующий аналог результата из~\cite {BF1}.

\begin {Sta}
Классы эквивалентности функций Морса $f\in F=F_{p,q,r}$ с
фиксированными кри\-ти\-чес\-кими точками на замкнутой ориентируемой
поверхности находятся во взаимно однозначном соответствии с
оснащенными нумерованными молекулами таких функций. В частности,
имеется лишь конечное число классов эквивалентности функций Морса
$f\in F=F_{p,q,r}$.
 \end {Sta}

\medskip
{\bf 5. Полиэдральные комплексы функций Морса и их разветвленные
накрытия.}
\begin {Def} \label {def:pol}
(А) Клеточный комплекс $X$ назовем {\it (кусочно евклидовым)
полиэдральным комп\-лексом}~\cite{BrHae}, если каждая его замкнутая
клетка $\bar\sigma$ снабжена метрикой, со\-гла\-со\-ванной с
индуцированной топо\-ло\-гией на $\bar\sigma$, и изометрична
некоторому выпуклому многограннику $P_{\sigma}$, причем изометрия
$\bar\sigma\to P_{\sigma}$ изометрично переводит все замкнутые клетки
$\bar\tau\subset\d\bar\sigma$ в грани многогранника $P_{\sigma}$.

(Б) Отображение $r\col\K\to K$ полиэдральных комплексов назовем {\it
правильным}, если его огра\-ни\-че\-ние на любую клетку
$\tilde\sigma$ комплекса $\K$ является изометрией на некоторую клетку
$\sigma$ комплекса $K$. Клетку $\tilde\sigma$ назовем {\it поднятием}
клетки $\sigma$ при $r$. В частности, $r$ является клеточным
отображением. Правильные биекции $K\to K$ назовем {\it
автоморфизмами} полиэдрального комплекса $K$.

(В) Пусть $\sigma,\tau\subset X$ -- два непересекающихся подмножества
топологического пространства $X$ (например, две открытые клетки
клеточного комплекса). Будем говорить, что $\sigma$ {\it примыкает к}
$\tau$ и писать $\tau\prec\sigma$ (и $\bar\tau\prec\bar\sigma$), если
$\tau\subset\d\sigma:=\bar\sigma\setminus\sigma$. Пишем
$\tau\preceq\sigma$, если $\tau\prec\sigma$ или $\tau=\sigma$.
 \end {Def}

\begin {Def} \label {def:cov}
Отображение $r\col\K\to K$ полиэдральных комплексов назовем {\it
раз\-ветв\-ленным накрытием}, если оно правильное (см.\
определение~\ref {def:pol}(Б)) и для любой клетки
$\tilde\tau\subset\K$ любая клетка $\sigma\subset K$, примыкающая к
клетке $\tau:=r(\tilde\tau)$ (см.\ определение~\ref {def:pol}(В)),
имеет поднятие $\tilde\sigma\subset\K$ (см.\ определение~\ref
{def:pol}(Б)), примыкающее к клетке $\tilde\tau$.
 \end {Def}

\begin {Thm} \label {thm:1}
Пусть количество седловых критических точек $q\ge1$. Существуют
$(q-1)$-мерный выпуклый многогранник $\P^{q-1}$ и $(q-1)$-мерные
поли\-эд\-ральные комплексы $\K=\K_{p,q,r}$ и $K=K_{p,q,r}$
(зависящие от чисел $p,r,q$ критических точек локальных минимумов,
максимумов и седловых точек), ассоци\-и\-ро\-ванные с пространством
$F=F_{p,q,r}$ функций Морса, и разветвленные накрытия $\K
\stackrel{r}{\longrightarrow} K \stackrel{r_0}{\longrightarrow}
\P^{q-1}$, такие что комплекс $K$ конечен и связен, и выполнены
следующие условия:

{\rm (А)} Пространство $F$ гомотопически эквивалентно полиэдральному
комплексу $\tilde K$.

{\rm (Б)} Клетки комплекса $\K$ (соответственно $K$) находятся во
взаимно однозначном соответ\-ствии с классами изотопности $[f]_\isot$
(соответственно классами эквивалентности $[f]$) функций Морса $f\in
F$. Размерность любой клетки равна $q-s(f)$, где $s(f)$ равно
количеству седловых критических значений функции $f\in F$, отвечающей
данной клетке. Две клетки $\tau,\sigma$ комплекса $\tilde K$
(соот\-вет\-ст\-венно $K$) примыкают друг к другу: $\tau\prec\sigma$
тогда и только тогда, когда соответ\-ст\-вующие им классы функций
Морса $[f]_\isot\leftrightarrow\sigma$,
$[g]_\isot\leftrightarrow\tau$ примыкают друг к другу как
подмно\-жества $F$ в $C^\infty$-топологии: $[f]_\isot\prec[g]_\isot$
(соответственно $[f]\prec[g]$, где $[f]\leftrightarrow\sigma$ и
$[g]\leftrightarrow\tau$).

{\rm (В)} Имеется правое действие группы $\G/\G^0$ на комплексе
$\tilde K$ автоморфизмами поли\-эдраль\-ного комплекса, согласованное
с естественным правым действием группы $\G$ на пространстве $F$.
Разветвленное накрытие $r\col\K\to K$ является $\G/\G^0$-инвариантным
и переводит друг в друга клетки $\tilde\sigma\to\sigma$, отвечающие
классам $[f]_\isot\subset[f]$ одной и той же функции Морса $f\in F$.
 \end {Thm}

Пункт (А) теоремы~\ref {thm:1} и утверждение о том, что $r_0$ --
разветвленное накрытие, не будут исполь\-зо\-ваны в настоящей работе;
их доказательство будет дано в следующих публикациях на основе~\cite
{KP1}.

\PProof{пунктов (Б,В) теоремы~\ref {thm:1}} Шаг 1. Опишем построение
выпуклого многогранника $\P^{q-1}$. Пусть $\P^{q-1}\subset\RR^q$ --
выпуклая оболочка множества точек $P_\pi:=\sum_{k=1}^q
\left(k-\frac{q+1}{2}\right)e_{\pi_k}$, $\pi\in\Sigma_q$, где
$e_1,\dots,e_q$ -- стандартный базис $\RR^q$. Известно~\cite{P}, что
$\P^{q-1}$ -- это $(q-1)$-мерный выпуклый многогранник в евклидовом
пространстве $\EE^{q-1}:=(e_1+\ldots+e_q)^\perp$, имеющий ровно $q!$
вершин $P_\pi$, $\pi\in\Sigma_q$, причем его $(q-s)$-мерные грани
находятся во взаимно однозначном соответствии с упорядоченными
разбиениями $J=(J_1,\dots,J_s)$ множества $\{1,\dots,q\}$ на $s$
непустых подмножеств $J_1,\dots,J_s$ (т.е.\ $\{1,\dots,q\}=J_1\sqcup
\ldots \sqcup J_s$), $1\le s\le q$. А именно, грань
$\tau^{q-s}_J\subset\P^{q-1}$, отвечающая разбиению $J$, -- это
выпуклая оболочка множества точек
$(\Sigma_{r_1}\times\Sigma_{r_2-r_1}\times\ldots\times\Sigma_{r_s-r_{s-1}})(P_\pi)$,
где числа $0=r_0<r_1<\ldots<r_{s-1}<r_s=q$ и перестановка
$\pi\in\Sigma_q$ однозначно определяются условиями
 \begin{equation}\label{eq:J}
J=(J_1,\dots,J_s), \quad J_1=\{\pi_1,\dots,\pi_{r_1}\}, \
J_2=\{\pi_{r_1+1},\dots,\pi_{r_2}\}, \ \ldots, \
J_s=\{\pi_{r_{s-1}+1},\dots,\pi_{r_s}\},
 \end{equation}
$\pi_1<\ldots<\pi_{r_1}$, $\pi_{r_1+1}<\ldots<\pi_{r_2}, \ldots,
\pi_{r_{s-1}+1}<\ldots<\pi_{r_s}$. Здесь
$\Sigma_{r_1}\times\Sigma_{r_2-r_1}\times\ldots\times\Sigma_{r_s-r_{s-1}}$
-- подгруппа группы $\Sigma_q$, отвечающая разбиению
$\{1,\dots,q\}=\{1,\dots,r_1\}\sqcup\{r_1+1,\dots,r_2\}\sqcup\ldots\sqcup\{r_{s-1}+1,\dots,r_s\}$,
и действие перестановки $\rho\in\Sigma_q$ на точке $P_\pi$ дает точку
$P_{\rho\pi}$, где $(\rho\pi)_i:=\pi_{\rho_i}$, $1\le i\le q$.

Если разбиение $\hat J$ получается из разбиения $J=(J_1,\dots,J_s)$
путем измельчения (т.е.\ разбиения некоторых множеств $J_k$ на
несколько подмножеств), будем писать $\hat J\prec J$. Из описания
граней многогранника $\P^{q-1}$ следует, что условие $\hat J\prec J$
равносильно $\tau_{\hat J}\prec\tau_J$ (см.\ определение~\ref
{def:pol}(В)).

Шаг 2. Для каждой функции Морса $f\in F$ рассмотрим набор $\bar
c=\bar c(f)=(c_1,\dots,c_q)\in\RR^q$ ее седловых критических значений
$c_i:=f(x_i)$, $1\le i\le q$. Сопоставим набору $\bar
c=(c_1,\dots,c_q)$ число $s(\bar c):=|\{c_1,\dots,c_q\}|$ различных
седловых значений и упорядоченное разбиение $J(\bar
c)=(J_1,\dots,J_s)$ множества $\{1,\dots,q\}$, определяемое
свойствами~(\ref {eq:J}) и $c_{\pi_1}=\ldots=c_{\pi_{r_1}} <
 c_{\pi_{r_1+1}}=\ldots=c_{\pi_{r_2}} < \ldots <
c_{\pi_{r_{s-1}+1}}=\ldots=c_{\pi_{r_s}}$. Сопоставим разбиению
$J(\bar c)$ и классу эквивалентности $[f]$ грань $\tau_{J(\bar
c)}\subset\P^{q-1}$.

Шаг 3. Покажем, что для любой функции $f\in F$ имеется биекция
$\delta[f]$ между множеством всех граней
$\tau'\prec\tau:=\tau_{J(\bar c(f))}$ и множеством всех классов
эквивалентности $[g]\succ[f]$ (см.\ определение~\ref {def:pol}(В)),
такая что $\delta{[f]}\col\tau'\mapsto[g]=:\delta_{\tau'}[f]$ при
$\tau'=\tau_{J(\bar c(g))}$. Это следует из следующих двух свойств:

1) для любого $\bar c\in\RR^q$ существует $\eps_0>0$, такое что (i)
для любого $\bar c'\in\RR^q$ со свойством $|\bar c'-\bar c|<\eps_0$
выполнено $J(\bar c')\preceq J(\bar c)$, и (ii) для любых
$\eps\in(0,\eps_0]$ и разбиения $\hat J\preceq J(\bar c)$ существует
$\bar c'\in\RR^q$ со свойствами $|\bar c'-\bar c|<\eps_0$ и $J(\bar
c')=\hat J$;

2) согласно~\cite[утверждение~1.1 и~\S3]{K}, любая функция $f\in F$
имеет окрестность $U$ в $F$, такую что для любых $g,g_1\in U$
равенства $[g]=[g_1]$ и $J(\bar c(g))=J(\bar c(g_1))$ равносильны.

Из этих свойств получаем, что из $[h]\succ[g]\succ[f]$ следует
$[h]\succ[f]$. Поэтому
 \begin{equation} \label {eq:*delta}
\delta_{\tau''}{[f]}=\delta_{\tau''}\delta_{\tau'}{[f]}
 \quad \mbox{для любых граней }\tau''\prec\tau'\prec\tau_{J(\bar c(f))}.
 \end{equation}

Шаг 4. Опишем построение полиэдрального комплекса $K$,
удовлетворяющего условиям пункта (Б), вместе с правильным
отображением $r_0\col K\to\P^{q-1}$. Рассмотрим метрическое
пространство $X:=\bigsqcup\limits_{[f]\in F/\sim} \ups_{[f]}$, где
$\ups_{[f]}$ является выпуклым многогранником, изометричным грани
${\tau_{J(\bar c(f))}}\subset\P^{q-1}$. Фиксируем отображение
$\pi\col X\to\P^{q-1}$, ограничение которого на каждый многогранник
$\ups_{[f]}$ является изометрией $\ups_{[f]}\to\tau_{J(\bar c(f))}$.
Очевидно, $\pi$ является правильным отобра\-же\-нием
поли\-эд\-ральных комплексов (см.\ определение~\ref {def:pol}(Б)).
Обозначим $\varphi_{[f]}:=(\pi|_{\ups_{[f]}})^{-1}\col\tau_{J(\bar
c(f))}\to\ups_{[f]}$.

Опишем (индукцией по $k\ge 0$) построение отношения эквивалентности
на множестве $X^{(k)}:=\bigsqcup\limits_{\dim\ups_{[f]}\le k}
\ups_{[f]}\subset X$ вместе с отображением $\pi_{k}\col
K^{(k)}\to(\P^{q-1})^{(k)}$, таких что
 \begin{equation} \label {eq:*k}
\pi_{k}\circ p_k=\pi|_{X^{(k)}}, \quad
p_{k}\circ\varphi_{[f]}|_{\tau'}=p_{k}\circ\varphi_{\delta_{\tau'}[f]}
\quad \mbox{для любых } f\in F,\ \dim\ups_{[f]}\le k,\quad \mbox{и }
\ \tau'\prec\tau_{J(\bar c(f))},
 \end{equation}
где $K^{(k)}$ -- множество классов эквивалентности в $X^{(k)}$,
$p_{k}\col X^{(k)}\to K^{(k)}$ -- каноническая проекция,
$\delta{[f]}$ -- биекция из шага~3. При $k=0$ различные точки считаем
не эквивалентными, определим $\pi_{0}$ формулой
$\pi_{0}(\ups_{[f]}):={\tau_{J(\bar c(f))}}$ при $\dim\ups_{[f]}=0$,
тогда выполнено~(\ref {eq:*k}) для $k=0$. Пусть $k\ge1$ и отношение
эквивалентности на $X^{(k-1)}$ с отображением $\pi_{k-1}$ уже
построены, причем $K^{(k-1)}$ является $(k-1)$-мерным полиэдральным
комплексом, $\pi_{k-1}$ -- правильным ото\-бра\-же\-нием и
вы\-пол\-нено~(\ref {eq:*k}) для $k-1$. Из~(\ref {eq:*delta}) и~(\ref
{eq:*k}) для $k-1$ следует, что для каждого $[f]$,
$\dim\ups_{[f]}=k$, имеется правильное вложение
$\varphi'_{[f]}\col\d{\tau_{J(\bar c(f))}}\to K^{(k-1)}$, такое что
$\varphi'_{[f]}|_{\tau'}=p_{k-1}\circ\varphi_{\delta_{\tau'}[f]}$ для
любого $\tau'\prec\tau_{J(\bar c(f))}$. Определим отношение
эквивалентности на $K^{(k-1)}\sqcup
\left(\bigsqcup\limits_{\dim\ups_{[f]}=k} \ups_{[f]}\right)$,
отождествляя каждую точку из $\d\ups_{[f]}$ с ее образом при
правильном вложении $\varphi'_{[f]}\circ\pi$. Тогда выполнено~(\ref
{eq:*k}), откуда $K^{(k)}$ -- $k$-мерный полиэдральный комплекс и
$\pi_k\col K^{(k)}\to(\P^{q-1})^{(k)}$ -- правильное отображение.

Таким образом, мы построили отношение эквивалентности $\sim_\glue$ на
всем $X=X^{(q-1)}$, поли\-эд\-раль\-ный комплекс
$K=K^{(q-1)}=X/\sim_\glue$ и правильное отображение
$r_0=\pi_{q-1}\col K\to\P^{q-1}$.

Шаг 5. Из утверждения и теоремы~\ref {thm:1}(Б) следует, что
полиэдральный комплекс $K$ конечен. Из результата о приведении
функций Морса к нормальной форме~\cite{Kmsb} следует, что $K$ связен.
Аналогично шагам 2--4 строится полиэдральный комплекс $\K$,
удовлетворяющий условиям пункта (Б), вместе с правильным отображением
$\K\to\P^{q-1}$ (для этого надо всюду в шагах 2--4 заменить
$[f],\ups_{[f]},X,X^{(k)}$,
$K^{(k)},\pi,\pi_k,p_k,\varphi_{[f]},\varphi'_{[f]}$ на
$[f]_\isot,\tilde\ups_{[f]_\isot},\tilde X,\tilde
X^{(k)},\K^{(k)},\tilde\pi,\tilde\pi_k,\tilde
p_k,\tilde\varphi_{[f]_\isot},\tilde\varphi'_{[f]_\isot}$).
Рассмот\-рим правое дей\-ствие группы $\G/\G^0$ на $\tilde
X:=\bigsqcup\limits_{[f]_\isot\in
F/\sim_\isot}\tilde\ups_{[f]_\isot}$, где элемент $h\G^0\in\G/\G^0$
действует по правилу
$h\G^0|_{\tilde\ups_{[f]_\isot}}:=\tilde\varphi_{[fh]_\isot}\circ\tilde\pi|_{\tilde\ups_{[f]_\isot}}\col\tilde\ups_{[f]_\isot}\to\tilde\ups_{[fh]_\isot}$;
тогда $X\approx\tilde X/(\G/\G^0)$. Это действие индуцирует действие
группы $\G/\G^0$ на $\K$ автоморфизмами полиэдрального комплекса (так
как отображения $\tilde\pi\col\tilde X\to\P^{q-1}$ и
$\delta_{\tau'}$, а потому и отношение эквивалентности $\sim_\glue$
на $\tilde X$, $\G/\G^0$-инвариантны). Поэтому композиция
правильного $\G/\G^0$-инвариантного отображения $\tilde X\to X$ и
правильного отобра\-же\-ния $X\to K$ индуцирует
правильное $\G/\G^0$-инвариантное отображение $r\col\K\to K$, такое
что $r(\tilde\sigma)=\sigma\leftrightarrow[f]$ при
$\tilde\sigma\leftrightarrow[f]_\isot$. Отсюда $r$ -- разветвленное
накрытие (см.\ определение~\ref {def:cov}).
 \qed

\begin{Not} \label {not:stab}
Для любой клетки $\hat\tau$ комплекса $\tilde K$ обозначим через
$\G^{\hat\tau}$ множество элементов $h\in\G/\G^0$, таких что
$\hat\tau h=\hat\tau$ (см.\ теорему~\ref {thm:1}(В)). Пусть $K^{(r)}$
-- $r$-мерный остов комплекса $K$.
\end{Not}

\begin{Thm} \label {thm:K}
Пусть $q\ge1$ и $f\in F$. Имеется эпиморфизм
$\mu\col\pi_1(K)\to\G_f/\D_f$. В частности, группа $\G_f/\D_f$ имеет
набор образующих $\mu([\gamma_1]),\dots,\mu([\gamma_\ell])$, где
$[\gamma_1],\dots,[\gamma_\ell]$ -- образующие $\pi_1(K)$.
\end{Thm}

\Proof Пусть $\tau\subset K$ и $\tilde\tau\subset\tilde K$ -- клетки
комплексов $K$ и $\tilde K$, отвечающие классам $[f]$ и $[f]_\isot$
(см.\ теорему~\ref {thm:1}(Б)). Без ограничения общности считаем, что
эти клетки нульмерны. Пусть $\K_f$ -- связная компонента комплекса
$\K$, содержащая клетку $\tilde\tau$. Рассмотрим правое действие
группы $\G_f/\G^0$ и ее подгруппы $\D_f/\G^0$ на $\K_f$ (см.\
теорему~\ref {thm:1}(В)). Так как разветвленное накрытие $\K_f\to K$
является $\G_f/\G^0$-инвариантным (см.\ там же), то
$K'_f:=\K_f/(\D_f/\G^0)$ -- полиэдральный комплекс, а проекция
$r'_f\col K'_f\to K\approx K'_f/(\G_f/\D_f)$ -- разветвленное
накрытие. В действи\-тель\-ности, $r'_f$ является накрытием, так как
$K'_f$ связен и действие на нем группы $\G_f/\D_f$ свободно (в силу
$\G^{\hat\tau}\subset\D_f/\G^0$). Поэтому имеется естественный
эпиморфизм $\mu\col\pi_1(K,\tau)\to\G_f/\D_f$, переводящий
гомотопический класс любой петли $\gamma\col[0,1]\to K$,
$\gamma(0)=\gamma(1)=\tau$, в элемент $h_\gamma\in\G_f/\D_f$, такой
что $\tilde\gamma(1)=\tilde\gamma(0)h_\gamma^{-1}$. Здесь
$\tilde\gamma\col[0,1]\to K'_f$ -- такое поднятие пути $\gamma$, что
$\tilde\gamma(0)=\tilde\tau(\D_f/\G^0)$.
 \qed

Опишем теперь образующие группы $\G_f/\G^0$ в терминах конечного
связного графа $K^{(1)}$. Пусть $T\subset K^{(1)}$ -- остовное дерево
графа $K^{(1)}$, пусть $\sigma_1,\dots,\sigma_n$ -- все ребра из
$K^{(1)}\setminus T$. Пусть $\tau_1,\dots,\tau_V$ и
$\sigma_1,\dots,\sigma_E$ -- все вершины и все ребра графа $K^{(1)}$
(каждое ребро снабдим произвольной ориентацией). Имеем $n=E-V+1$.
Пусть $S\col T\to\tilde K$ -- любое непрерывное поднятие дерева $T$,
такое что $S(\tau)=\tilde\tau$ (здесь $\tau,\tilde\tau$ как в
доказательстве теоремы~\ref {thm:K}), и пусть $\hat\sigma_e$ -- такое
поднятие ребра $\sigma_e$, что $\hat\sigma_e(0)=S(\sigma_e(0))$,
$1\le e\le n$. Имеем $\hat\sigma_e(1)=S(\sigma_e(1))h_e$ для
некоторого $h_e\in\G_f/\G^0$, $1\le e\le n$. Элементы
$h_1,\dots,h_n\in\G_f/\G^0$ назовем {\it $T$-дополнительными
элементами}.

\begin{Thm} [М.~Басманова и Е.~Кудрявцева, 1999] \label {thm:BK}
Пусть $q\ge1$ и $f\in F$. Группа $\G_f/\G^0$ имеет конечную систему
образующих $A_1\cup\ldots\cup A_V\cup\{h_1,\dots,h_n\}$, где $A_v$ --
конечная система образующих группы $\G^{S(\tau_v)}$, $1\le v\le V$
{\rm (см.\ обозначение~\ref {not:stab})}, $h_1,\dots,h_n\in\G_f/\G^0$
-- $T$-дополнительные элементы. Для минимального числа образующих
верно $\rank(\G_f/\G^0)\le (q+g-1)V+n=(q+g-2)V+E+1$, где $V$ и $E$ --
количества вершин и ребер графа $K^{(1)}$, $n=E-V+1$, $g$ -- род
поверхности $M$.
\end{Thm}

\Proof Пусть $h\in\G_f/\G^0$. Тогда, в обозначениях доказательства
теоремы~\ref {thm:K}, существуют петля $\gamma\col[0,1]\to K^{(1)}$ и
ее поднятие $\tilde\gamma$ в $\tilde K$, такие что
$\tilde\gamma(0)=\tilde\tau$ и $\tilde\gamma(1)=\tilde\tau h$. Пусть
$\gamma=\sigma_{e_1}^{\eps_1}\cdot\ldots\cdot\sigma_{e_N}^{\eps_N}$
-- разложение петли $\gamma$ в произведение ориентированных ребер
комплекса $K$, где $\eps_i\in\{1,-1\}$ и $e_i\in\{1,\dots,E\}$, $1\le
i\le N$, и пусть
$\tilde\gamma=\tilde\sigma_1^{\eps_1}\cdot\ldots\cdot\tilde\sigma_N^{\eps_N}$
-- соответствующее разложение. Обозначим
$\tilde\tau_i:=\tilde\sigma_i^{\eps_i}(1)$ при $1\le i\le N$;
$\hat\sigma_e:=S(\sigma_e)$ и $h_e:=1\in\G_f/\G^0$ при $n<e\le E$.
Имеем $\sigma_{e_i}^{\eps_i}(1)=\tau_{v_i}$ для некоторого
$v_i\in\{1,\dots,V\}$;
$\tilde\sigma_i^{\eps_i}=\hat\sigma_{e_i}^{\eps_i}\tilde h_i$ для
некоторого $\tilde h_i\in\G/\G^0$ ($1\le i\le N$).

Из
$\tilde\tau=\tilde\gamma(0)=\tilde\sigma_1^{\eps_1}(0)=\hat\sigma_{e_1}^{\eps_1}(0)\tilde
h_1$ имеем $\tilde\tau=\tilde\tau h_{e_1}^{\frac{1-\eps_1}2}\tilde
h_1$, откуда $\tilde h_1\in
h_{e_1}^{\frac{\eps_1-1}2}\G^{\tilde\tau}$ (см.\ обозначение~\ref
{not:stab}). При $1\le i<N$ из
$\tilde\sigma_{i}^{\eps_i}(1)=\hat\sigma_{e_i}^{\eps_i}(1)\tilde
h_{i}$ и
$\tilde\sigma_{i+1}^{\eps_{i+1}}(0)=\hat\sigma_{e_{i+1}}^{\eps_{i+1}}(0)\tilde
h_{i+1}$ имеем
$\tilde\tau_{i}=S(\tau_{v_i})h_{e_i}^{\frac{\eps_i+1}2}\tilde h_{i}$
и
$\tilde\tau_{i}=S(\tau_{v_i})h_{e_{i+1}}^{\frac{1-\eps_{i+1}}2}\tilde
h_{i+1}$, откуда $\tilde h_{i+1}\tilde h_i^{-1}\in
h_{e_{i+1}}^{\frac{\eps_{i+1}-1}2}\G^{S(\tau_{v_i})}h_{e_i}^{\frac{\eps_i+1}2}$.
Из $\tilde\sigma_{N}^{\eps_N}(1)=\hat\sigma_{e_N}^{\eps_N}(1)\tilde
h_{N}$ и $\tilde\gamma(1)=\tilde\tau h$ имеем
$\tilde\tau_N=S(\tau_{v_N})h_{e_N}^{\frac{\eps_N+1}2}\tilde
h_N=\tilde\tau h_{e_N}^{\frac{\eps_N+1}2}\tilde h_N$ и
$\tilde\tau_N=\tilde\tau h$, откуда $h\tilde
h_N^{-1}\in\G^{\tilde\tau}h_{e_N}^{\frac{\eps_N+1}2}$. Поэтому
$$
 h\ =\ (h\tilde h_N^{-1})\ (\tilde h_N\tilde h_{N-1}^{-1})\ \ldots\ (\tilde h_2\tilde h_1^{-1})\ \tilde h_1
 \ \in \ \G^{\tilde\tau}\ h_{e_N}^{\eps_N}\ \G^{S(\tau_{v_{N-1}})}\ h_{e_{N-1}}^{\eps_{N-1}}\
 \ldots\ h_{e_2}^{\eps_2}\ \G^{S(\tau_{v_1})}\ h_{e_1}^{\eps_1}\ \G^{\tilde\tau},
$$
т.е.\ $h$ есть произведение степеней элементов из $A_1\cup\ldots\cup
A_V\cup\{h_1,\dots,h_n\}$. Оценка $\rank(A_v)\le q+g-1$ легко
доказывается, см.\ замечание перед следствием.
 \qed

\medskip
Автор приносит благодарность Д.M.\ Афанасьеву, М.\ Басмановой, Ю.М.\
Бурману, М.\ Концевичу, Д.А.\ Пермякову, Л.\ Фадеевой, А.Т.\ Фоменко
и Х.\ Цишангу за полезные замечания и обсуждения.

Работа частично поддержана грантом РФФИ \No~10--01--00748-а,
грантом про\-граммы ``Ведущие научные школы РФ'' НШ-3224.2010.1,
грантом программы ``Развитие научного потенциала высшей школы'' РНП
2.1.1.3704 <<Современная дифференциальная геометрия, топология и
приложения>> и грантом ФЦП <<Научные и научно-педагогические кадры
инновационной России>> (контракты \No~02.740.11.5213 и \No~14.740.11.0794).


\noindent
Mathematics and Mechanics Department of Moscow State University \\
{\it E-mail address: } eakudr@mech.math.msu.su

\end{document}